\documentclass{article}
\usepackage{amssymb,latexsym,amsmath,amsthm,mathrsfs}
\usepackage{graphicx}
\newcommand{\comment}[1]{}

\begin{document}
\title{Various considerations on hypergeometric series\footnote{Presented
to the St. Petersburg Academy on August 19, 1776.
Originally
published as
{\em Variae considerationes circa series hypergeometricas},
Nova Acta Academiae Scientarum Imperialis Petropolitinae \textbf{8} (1792),
3--14.
E661 in the Enestr{\"o}m index.
Translated from the Latin by Jordan Bell,
Department of Mathematics, University of Toronto, Toronto, Ontario, Canada.
Email: jordan.bell@gmail.com}}
\author{Leonhard Euler}
\date{}
\maketitle

\S 1. Given this product continuing to infinity:
\[
\frac{a(a+2b)}{(a+b)(a+b)}\cdot \frac{(a+2b)(a+4b)}{(a+3b)(a+3b)}
\cdot \frac{(a+4b)(a+6b)}{(a+5b)(a+5b)}
\cdot \frac{(a+6b)(a+8b)}{(a+7b)(a+7b)} \textrm{\& c.}=\frac{P}{Q},
\]
it is known that\footnote{Translator: $P=\frac{1}{2b}B(\frac{a+b}{2b},\frac{1}{2})$ and $Q=\frac{1}{2b}B(\frac{a}{2b},\frac{1}{2})$, where $B(x,y)$
is the Beta function. $B(x,q)=\frac{\Gamma(x)\Gamma(y)}{\Gamma(x+y)}$,
where $\Gamma$ is the Gamma function.
See Chapter XII of Whittaker and Watson. In this paper Euler does not use $\Gamma$ to represent the usual Gamma function.}
\[
P=\int \frac{x^{a+b-1}\partial x}{\surd(1-x^{2b})} \quad \textrm{and} \quad
Q=\int \frac{x^{a-1}\partial x}{\surd(1-x^{2b})}
\]
with the integrals extended from $x=0$ to $x=1$;
it is worthwhile to note here that the member of this product corresponding
to the index $i$ is
\[
\frac{(a+(2i-2)b)(a+2ib)}{(a+(2i-1)b)(a+(2i-1)b)}.
\]

\S 2. Now for the study of this product, let us consider the following
indefinite product, in which the number of factors $=n$, and let us put  
\[
a(a+2b)(a+4b)(a+6b)\cdots(a+(2n-2)b)=\Delta:n,
\]
if indeed this product, with $a$ and $b$ given number,
is seen as a certain function of $n$;
then from its nature it is clear that
\[
\Delta:(n+1)=\Delta:n\cdot(a+2nb),
\]
similarly
\[
\Delta:(n+2)=\Delta:(n+1)\cdot (a+(2n+2)b),
\]
and so on. Then if $i$ denotes an infinitely large number, it will be
\[
\Delta:i=a(a+2b)(a+4b)(a+6b)\cdots(a+(2i-2)b)
\]
whence it even follows
\begin{eqnarray*}
\Delta:(i+1)&=&\Delta:i\cdot(a+2ib),\\
\Delta:(i+2)&=&\Delta:i\cdot(a+2ib)(a+(2i+2)b),\\
\Delta:(i+3)&=&\Delta:i\cdot(a+2ib)(a+(2i+2)b)(a+(2i+4)b)\\
&&\textrm{etc.},
\end{eqnarray*}
where the factors that follow can be considered as equal to each other;
thus it can in general be put
\[
\Delta:(i+n)=\Delta:i\cdot(a+2ib)^n,
\]
where,
since the factor $(a+2ib)$ is approximately the next, 
by this reason each of the following have been obtained,
from which we can also put in general
\[
\Delta:(i+n)=(\alpha+2ib)^n \cdot \Delta:i,
\]
where $\alpha$ denotes any finite number, 
which of course vanishes next to $2ib$.

\S 3. Let us now work out the computation for the case of the indefinite
product in which $n=\frac{1}{2}$, and let us call $\Delta:\frac{1}{2}=k$,
whose value can indeed always be assigned approximately by means of the method
of interpolation. Here, from the above it will therefore be
\begin{eqnarray*}
\Delta:\Big(1+\frac{1}{2}\Big)&=&k(a+b),\\
\Delta:\Big(2+\frac{1}{2}\Big)&=&k(a+b)(a+3b),\\
\Delta:\Big(3+\frac{1}{2}\Big)&=&k(a+b)(a+3b)(a+5b)\\
&&\textrm{etc.},
\end{eqnarray*}
whence by progressing to infinity it will be
\[
\Delta:\Big(i+\frac{1}{2}\Big)=k(a+b)(a+3b)(a+5b)\cdots(a+(2i-1)b).
\]

\S 4. Since we already gave above the formula for $\Delta:(i+n)$,
by now putting $n=\frac{1}{2}$ we will also have
\[
\Delta:\Big(i+\frac{1}{2}\Big)=\Delta:i \surd(a+2ib),
\]
and thus for the one formula $\Delta:\Big(i+\frac{1}{2}\Big)$ we have found
two different expressions; from them one constructs this equation
\[
\Delta:i\surd(\alpha+2ib)=k(a+b)(a+3b)(a+5b)\cdots(a+(2i-1)b),
\]
and we can thus conclude that the value of this infinite product
is
\[
(a+b)(a+3b)(a+5b)\cdots(a+(2i-1)b)=\frac{\Delta:i\surd(\alpha+2ib)}{k},
\]
and thus a relation is revealed between this product and that which
we expressed above by $\Delta:i$. 
It is also appropriate to note here that
the factors of this product are the same as those which constitute
the denominator of the product first proposed;\footnote{Translator: That is
the factors here are the same as the factors of $Q$, though here they 
are simple and are squared in $Q$.} therefore we will be able to express
the numerator and denominator of this product respectively by the values
just found
\[
\Delta:i \quad \textrm{and} \quad \frac{\Delta:i\surd(\alpha+2ib)}{k}
\]

\S 5. Furthermore, the numerator of the proposed product expanded
to infinity can be thus represented
\[
a(a+2b)^2(a+4b)^2\cdots (a+(2i-2)b)^2(a+2ib),
\]
where the first and last factors are alone, while all the others are squares.
Since it will thus be
\[
(\Delta:i)^2=(a)^2(a+2b)^2(a+4b)^2(a+6b)^2\cdots (a+(2i-2)b)^2,
\]
it is evident that the numerator is $\frac{(\Delta:i)^2}{a}(a+2ib)$. For the
denominator, it is also clear that it is equal to the square of the other
product $(a+b)(a+3b)$ etc.; since the value of this has been found to be
\[
\frac{\Delta:i\surd(\alpha+2ib)}{k},
\]
the denominator will be
\[
\frac{(\Delta:i)^2(\alpha+2ib)}{kk};
\] 
then with these values substituted into the fraction $\frac{P}{Q}$ shown above,
we have obtained this equation
\[
\frac{P}{Q}=\frac{\frac{(\Delta:i)^2(a+2ib)}{a}}{\frac{(\Delta:i)^2(\alpha+2ib)}{kk}}
=\frac{kk(a+2ib)}{a(\alpha+2ib)}=\frac{kk}{a}.
\]
Therefore from this equation the true value of the interpolated
formula $k=\Delta:\frac{1}{2}$ can be found at once, since it will be
\[
\Delta:\frac{1}{2}=\surd\frac{aP}{Q}
\]
and then in turn the following
\begin{eqnarray*}
\Delta:\Big(1+\frac{1}{2}\Big)&=&(a+b)\surd\frac{aP}{Q},\\
\Delta:\Big(2+\frac{1}{2}\Big)&=&(a+b)(a+3b)\surd\frac{aP}{Q},\\
\Delta:\Big(3+\frac{1}{2}\Big)&=&(a+b)(a+3b)(a+5b)\surd\frac{aP}{Q}\\
&&\textrm{etc.},
\end{eqnarray*}
and this interpolation is itself all the more noteworthy, since it at once
provides the true value of the interpolated terms without approximation.

\S 6. But let us further consider that infinite product in which both
types of factors are joined, and let us put
\[
a(a+b)(a+2b)(a+3b)\cdots(a+(i-1)b)=\Gamma:i.
\]
It will be
\[
\Gamma:2i=a(a+b)(a+2b)(a+3b)\cdots(a+(2i-1)b),
\]
which is clearly the product of the above two, so it will be
\[
\Gamma:2i=\frac{(\Gamma:i)^2\surd(\alpha+2ib)}{k};
\]
thus if we want to make use of the form $\Gamma:2i$,
we can assign the values of both the preceding from it, so that it will be
\[
\Delta:i=\surd\frac{k\cdot \Gamma:2i}{\surd(\alpha+2ib)},
\]
which is the value of the first product
\[
a(a+2b)(a+4b)(a+6b)\, \textrm{etc.};
\]
and indeed the value of the latter product
\[
(a+b)(a+3b)(a+5b) \, \textrm{etc.}
\]
will be
\[
\surd\frac{\Gamma:2i\surd(\alpha+2ib)}{k}.
\]

\S 7. So far we have contemplated three products extending to infinity
and their relations to each other. To help us accurately examine these, let us
repeat them again for viewing
\begin{eqnarray*}
\textrm{I.}&&a(a+b)(a+2b)(a+3b)\cdots(a+(i-1)b)=\Gamma:i,\\
\textrm{II.}&&a(a+2b)(a+4b)(a+6b)\cdots(a+(2i-2)b)=\Delta:i,\\
\textrm{III.}&&(a+b)(a+3b)(a+5b)\cdots(a+(2i-1)b)=\Theta:i,
\end{eqnarray*}
and we have already found that\footnote{Translator: In \S 4.}
\[
\Theta:i=\frac{\Delta:i\surd(\alpha+2ib)}{k};
\]
next indeed we can express both $\Delta:i$ and $\Theta:i$ in the following
way with the function $\Gamma:2i$,
\[
\Delta:i=\surd\frac{k\cdot\Gamma:2i}{\surd(\alpha+2ib)} \quad \textrm{and} \quad
\Theta:i=\surd\frac{\Gamma:2i \surd(\alpha+2ib)}{k},
\]
from which it is clear that
\[
\Gamma:2i=\Delta:i \cdot \Theta:i;
\]
where it should be remembered that $k=\Delta:\frac{1}{2}$,
which namely ought to be defined from the second form by considering the
series\footnote{Translator: It looks like Euler is just writing down again the series
that we want to find the interpolated value $k$ from.}
\[
a, \, a(a+2b), \, a(a+2b)(a+4b), \, a(a+2b)(a+4b)(a+6b) \quad \textrm{etc.},
\]

\S 8. Now let us apply to these forms a general method for summing progressions
of all types using their general term,
which is such that given some series $A,B,C,D,E$ etc., whose term 
corresponding to the indefinite index $x$ is $=X$, with sum 
\[
A+B+C+D+\cdots+X,
\]
which we shall call $=S$, it will be\footnote{Translator: This is the Euler-Maclaurin summation formula.}
\[
S=\int X\partial x+\frac{1}{2}X+\frac{1}{1\cdot 2\cdot 3}\cdot \frac{1}{2}\frac{\partial X}{\partial x}
-\frac{1}{1\cdot 2\cdots 5}\cdot\frac{1}{6}\frac{\partial^3 X}{\partial x^3}
+\frac{1}{1\cdot 2\cdots 7}\cdot \frac{1}{6} \frac{\partial^5 X}{\partial x^5}
-\textrm{etc.},
\]
where the fractions $\frac{1}{2},\frac{1}{6},\frac{1}{6},\frac{3}{10},\frac{5}{6}$ etc. are the Bernoulli numbers.

\begin{center}
{\Large Development of the first form}
\[
a(a+b)(a+2b)(a+3b)\cdots(a+(i-1)b)=\Gamma:i.
\]
\end{center}

\S 9. Since here the number of factors is considered infinite,
to use this summation method we consider the same form with a number of terms
$=x$, a finite constant, and in this way let us set
\[
a(a+b)(a+2b)(a+3b)\cdots(a+(x-1)b)=\Gamma:x.
\]
Now indeed, to get a product in place of this series, let us take the logarithms, and it will be
\[
l\Gamma:x=la+l(a+b)+l(a+2b)+l(a+3b)+\cdots+l(a+(x-1)b);
\]
whose sum shall therefore be explored, which gives the logarithm 
of the formula $\Gamma:x$ and hence the formula $\Gamma:x$
itself.\footnote{Translator: That is, by knowing the logarithm of function we
know the function.}
If one then sets $x=i$ in this, the formula $\Gamma:i$ will be obtained,
whose value has been our main interest in the above. Therefore
comparing this with the series stated in the most general
form,\footnote{Translator: The series $S$ in the Euler-Maclaurin
summation formula} it will be
\[
X=l(a+(x-1)b)
\]
and then the sum
\[
S=l\Gamma:x,
\]
or it will be
\[
X=l(a-b+bx),
\]
from which one get
\[
\int X\partial x=\int \partial xl(a-b+bx).
\]

\S 10. Therefore since it is
\[
\int \partial zlz=zlz-z
\]
and so
\[
\int \partial yl(a+y)=(a+y)l(a+y)-(a+y),
\]
writing now $bx$ in place of $y$ it will be
\[
\int b\partial xl(a+bx)=(a+bx)l(a+bx)-a-bx
\]
and hence
\[
\int \partial xl(a+bx)=\frac{(a+bx)}{b}l(a+bx)-\frac{a}{b}-x,
\]
from which we get that in our case it is
\[
\int X\partial x=\frac{(a-b+bx)}{b}l(a-b+bx)-\frac{a}{b}+1-x,
\]
where in the last part the constant $\frac{a}{b}-1$ can be omitted,
since the expression requires for itself an indefinite constant quantity,
which ought to be defined later from the nature of the series. Next it will indeed be
\[
\frac{\partial X}{\partial x}=\frac{b}{a-b+bx},
\]
then indeed in turn
\[
\frac{\partial^3 X}{\partial x^3}=\frac{2b^3}{(a-b+bx)^3}, \quad
\frac{\partial^5 X}{\partial x^5}=\frac{2\cdot 3\cdot 4b^5}{(a-b+bx)^5} \quad
\textrm{etc.}
\]
By calling these values to use it will be
\begin{eqnarray*}
l\Gamma:x&=&A+\Big(\frac{a}{b}-\frac{1}{2}+x\Big)l(a-b+bx)-x+\frac{1}{1\cdot 2\cdot 3}\cdot \frac{1}{2}\cdot \frac{b}{a-b+bx}\\
&&-\frac{1}{3\cdot 4\cdot 5}\cdot \frac{1}{6}\cdot \frac{b^3}{(a-b+bx)^3}+
\frac{1}{5\cdot 6\cdot 7}\cdot \frac{1}{6}\cdot \frac{b^5}{(a-b+bx)^5}\\
&&-\frac{1}{7\cdot 8\cdot 9}\cdot \frac{3}{10}\cdot \frac{b^7}{(a-b+bx)^7}
+\frac{1}{9\cdot 10\cdot 11}\cdot \frac{5}{6}\cdot \frac{b^9}{(a-b+bx)^9}
-\textrm{etc.},
\end{eqnarray*}
where the letter $A$ denotes a constant to be defined from the nature of
the series.

\S 11. Moreover, this constant $A$ should be determined from a case in
which the sum of the series is known, which can thus be done in the case $x=0$,
because of course 
the sum ought to turn out to be equal to
nothing;\footnote{Translator: Since $\Gamma:0=1$} therefore it will be
\begin{eqnarray*}
-A&=&\Big(\frac{a}{b}-\frac{1}{2}\Big)l(a-b)+\frac{1}{1\cdot 2\cdot 3}\cdot \frac{1}{2}
\cdot \frac{b}{a-b}-\frac{1}{3\cdot 4\cdot 5}\cdot \frac{1}{6}\cdot \frac{b^3}{(a-b)^3}\\
&&+\frac{1}{5\cdot 6\cdot 7}\cdot \frac{1}{6}\cdot \frac{b^5}{(a-b)^5}-\textrm{etc.}
\end{eqnarray*}
Seeing that this series is hardly convergent, and in the case $b=a$ all the terms
become infinite, then apparently nothing can be done with it. But if
on the other hand we were to take $x=1$, this sum would turn out to
be $=la$,\footnote{Translator: Since $\Gamma:1=a$.}
from which it is likewise difficult to conclude anything for our purpose,
for this always leads to an infinite series whose sum must eventually be
explored. Perhaps indeed in this work some use could be made of series
involving the Bernoulli numbers, which I have previously commented on,\footnote{Translator: E393, ``De summis serierum numeros Bernoullianos involventium'' (1768).}
but I shall not further take up this question here.

\S 12. Since for our present investigation we are primarily interested in
the value $\Gamma:i$, it suffices to at once put an infinite number in place
of $x$. Therefore let $x=i$, $i$ denoting an infinitely large number, and
our equation takes this form
\[
l\Gamma:i=A+\Big(\frac{a}{b}-\frac{1}{2}+i\Big)l(a-b+bi)-i,
\]
where the constant $A$ is freely determined,
and here we can consider it as known. 
Therefore now returning to numbers,\footnote{Translator: We can consider
the above equation as a sum of logarithms. By taking the exponential of both sides
we get rid of the logarithms and have ``numbers'', not logarithms.}
in which we understand $lA$ as written in place of $A$, we arrive
at this expression
\[
\Gamma:i=A(a-b+bi)^{\frac{a}{b}-\frac{1}{2}+i}e^{-i}.
\]
It is convenient here to represent separately the power with the exponent
$i$ in this way
\[
\Gamma:i=A(a-b+bi)^{\frac{a}{b}-\frac{1}{2}}(a-b+bi)^i e^{-i}.
\]

\begin{center}
{\Large Development of the two other formulas}
\end{center}

\S 13. The second form differs from the first in that here $2b$ needs to be written
in place of $b$, from which we can avoid developing it from scratch; 
and indeed let us write $B$ in place of the constant
$A$, since it is not yet apparent how the letter $b$ is involved in the constant $A$. Thus we will now have
\[
\Delta:i=B(a-2b+2bi)^{\frac{a}{2b}-\frac{1}{2}}(a-2b+2bi)^ie^{-i}.
\]
In a similar way, it is evident that the third form arises from the second form
if now one writes $a+b$ in place of $a$, where by introducing the constant
$C$ in place of $B$ we will have
\[
\Theta:i=C(a-b+2bi)^{\frac{a}{2b}}(a-b+2bi)^i e^{-i}.
\]
Here it should be noted that the letter $e$ is taken here as that number
whose hyperbolic logarithm $=1$.

\begin{center}
{\Large Conclusions that arise from this}
\end{center}

\S 14. Let us now see how these new determinations
can be compared to those found with the relations found above; for with these
new values it will be
\[
\Gamma:2i=A(a-b+2bi)^{\frac{a}{b}-\frac{1}{2}}(a-b+2bi)^{2i}e^{-2i}.
\]
Since we have found that
\[
\Gamma:2i=\Delta:i \cdot \Theta:i,
\]
if we substitute the values found in this way everywhere, we will have
for the first equation the product
\begin{eqnarray*}
&&\Delta:i \cdot \Theta:i\\
&=&BC(a-2b+2bi)^{\frac{a}{2b}-\frac{1}{2}}(a-b+2bi)^{\frac{a}{2b}}(a-2b+2bi)^i
(a-b+2bi)^i e^{-2i};
\end{eqnarray*}
because this should be equal to the value of $\Gamma:2i$, if we divide
both sides by the common factors, we will obtain this equation
\[
A(a-b+2bi)^{\frac{a}{2b}-\frac{1}{2}}(a-b+2bi)^i=BC(a-2b+2bi)^{\frac{a}{2b}-\frac{1}{2}}(a-2b+2bi)^i.
\]

\S 15. Let us divide both sides of this equation by $(a-2b+2bi)^i$, and since
it will be
\[
\frac{a-b+2bi}{a-2b+2bi}=1+\frac{b}{a-2b+2bi}=1+\frac{1}{2i}
\]
where $i$ is an infinite number, by the ordinary resolution it will be
\[
\Big(1+\frac{1}{2i}\Big)^i=e^{\frac{1}{2}},
\]
with which our equation reduces to this form
\[
A(a-b+2bi)^{\frac{a}{2b}-\frac{1}{2}}\cdot e^{\frac{1}{2}}=BC(a-2b+2bi)^{\frac{a}{2b}-\frac{1}{2}},
\]
where the last factors cancel each other out, since
\[
\Big(\frac{a-b+2bi}{a-2b+2bi}\Big)^{\frac{a}{2b}-\frac{1}{2}}
=\Big(1+\frac{1}{2i}\Big)^{\frac{a}{2b}-\frac{1}{2}}=1,
\]
so that this simple equality is reached
\[
Ae^{\frac{1}{2}}=BC.
\]

\S 16. Next, since we have found above that
\[
\Theta:i=\frac{\Delta:i \surd(\alpha+2ib)}{k}
\]
or
\[
\frac{\Theta:i}{\Delta:i}=\frac{\surd(\alpha+2ib)}{k},
\]
let us divide the value found for $\Theta:i$ by $\Delta:i$, and we will get
\[
\frac{\Theta:i}{\Delta:i}=\frac{C}{B}\surd(a-2b+2bi)\Big(\frac{a-b+2bi}{a-2b+2bi}\Big)^i
=\frac{C}{B}\surd e(a-2b+2bi).
\]
Therefore it will be
\[
\frac{\surd(\alpha+2ib)}{k}=\frac{C}{B}\surd e(a-2b+2bi)
\]
or
\[
\frac{1}{k}=\frac{C}{B}\surd \frac{e(a-2b+2bi)}{\alpha+2ib}=\frac{C}{B}\surd e,
\]
or it will be
\[
B=Ck\surd e.
\]

\S 17. Therefore we have got two relations between the three constants $A,B,C$,
such that if one of them were known, from it the other two could be defined.
For since it is
\[
A=\frac{BC}{\surd e} \quad \textrm{and} \quad B=Ck\surd e,
\]
if we consider the constant $A$ as already known, the other two can be
determined in turn. For since $B=Ck\surd e$, this value subsituted into
the first equation gives $A=CCk$, from which one gets $C=\surd \frac{A}{k}$,
and then in turn $B=\surd k Ae$. It is still not clear however how the constant
$A$ can be absolutely determined, and thus it will return to the summation
of the logarithmic series which we indicated above with the letter $A$,\footnote{Translator: $A$ is defined in \S 11. I do not know what Euler means by ``logarithmic series''.}
where however $lA$ needs to be written in place of $A$. And likewise
we find
that if the two other forms are developed in the same way
into logarithmic series, the constants used in them, namely $lB$ and $lC$,
can be found at once.

\S 18. It still remains for us to add a bit to the value of the letter
$k$, which we already mentioned above must be found by
interpolation.\footnote{Translator: Here Euler seems to be showing a different argument
to find the value of $k$. Indeed he already found $k$ in \S 5. However I do
not see how this derivation is different from \S 5.}
Still however this letter can be determined absolutely
by certain quadratures, from the comparison of the formulas
$\Delta:i$ and $\Theta:i$. For since it is
\[
k=\frac{\Delta:i}{\Theta:i}\surd(\alpha+2ib)
\]
and hence
\[
kk=\frac{(\Delta:i)^2(\alpha+2ib)}{(\Theta:i)^2},
\]
if we substitute infinite products in place of $\Delta:i$ and $\Theta:i$,
since both are constructed from $i$ factors, the single factor $\alpha+2ib$ 
can be absorbed into the numerator,\footnote{Translator: It looks like Euler
is saying that since there are infinitely many infinitely large factors on the
numerator and the denominator, we can just ignore the single factor
$\alpha+2ib$.} and let us express the first factor of the numerator separately;
with this done we will be led to the following product to be determined:
\[
kk=a\cdot \frac{a(a+2b)(a+2b)(a+4b)(a+4b)(a+6b)}{(a+b)(a+b)(a+3b)(a+3b)(a+5b)(a+5b)}\cdot\, \textrm{etc.}
\]

\S 19. So that we can find true value of this infinite product, it should
be remembered that if the letters $P$ and $Q$ denote the following integral
fromulas
\[
P=\int \frac{x^{p-1}\partial x}{(1-x^n)^{1-\frac{m}{n}}} \quad \textrm{and} \quad
Q=\int \frac{x^{q-1}\partial x}{(1-x^n)^{1-\frac{m}{n}}},
\]
where the integrals are understood to extend from $x=0$ to $x=1$,
then we have this infinite product
\[
\frac{P}{Q}=\frac{q(m+p)}{p(m+q)}\cdot \frac{(q+n)(m+p+n)}{(p+n)(m+q+n)}
\cdot \frac{(q+2n)(m+p+2n)}{(p+2n)(m+q+2n)}\cdot \, \textrm{etc.},
\]
This product can easily be reduced to our form by taking
\[
q=a, \quad p=a+b, \quad m=b, \quad n=2b,
\]
so that in our case it becomes
\[
P=\int \frac{x^{a+b-1}\partial x}{\surd(1-x^{2b})} \quad \textrm{and}
\quad Q=\int \frac{x^{a-1}\partial x}{\surd(1-x^{2b})};
\]
it will then indeed be
\[
kk=\frac{aP}{Q}
\]
and so
\[
k=\surd\frac{aP}{Q},
\]
and thus we have elicited the same value $k$ in a different way than 
we produced it above.

\S 20. While it is $k=\Delta:\frac{1}{2}$, we can assign
values in a similar way to the two other forms $\Gamma:\frac{1}{2}$
and $\Theta:\frac{1}{2}$. For since the form $\Gamma$ arises from the form
$\Delta$ if one writes $\frac{1}{2}b$ in place of $b$,
and the form $\Theta$ from $\Delta$ if $a+b$ is written in place of
$a$, with this observed it will be
\[
\Gamma:\frac{1}{2}=\sqrt{a\frac{\int \frac{x^{a
+\frac{1}{2}b-1}\partial x}{\surd(1-x^b)}}{\int \frac{x^{a-1}\partial x}{\surd(1-x^b)}}}
\]
and
\[
\Theta:\frac{1}{2}=\sqrt{(a+b)\frac{\int \frac{x^{a+2b-1}\partial x}{\surd(1-x^{2b})}}{\int \frac{x^{a+b-1}\partial x}{\surd(1-x^{2b})}}}
\]
It is easy to see that the value of $\Theta:\frac{1}{2}$
can be introduced in our calculations like $\Delta:\frac{1}{2}$,
since it is
\[
\Delta:\frac{1}{2}\cdot \Theta:\frac{1}{2}=a.
\]
For on the other hand, taking these integral values 
yields
\[
\Delta:\frac{1}{2}\cdot \Theta:\frac{1}{2}=\sqrt{\frac{a(a+b)\int
\frac{x^{a+2b-1}\partial x}{\surd(1-x^{2b})}}{\int
\frac{x^{a-1}\partial x}{\surd(1-x^{2b})}}};
\]
and from a well known reduction of these integrals it turns out that
\[
\int \frac{x^{a+2b-1}\partial x}{\surd(1-x^{2b})}=\frac{a}{a+b}\int \frac{x^{a-1}\partial x}{\surd(1-x^{2b})},
\]
namely for the terms of integration $x=0$ and $x=1$, and thus it is clear that
\[
\Delta:\frac{1}{2}\cdot \Theta:\frac{1}{2}=a.
\]
However, there seems to be no way to
express the value of $\Gamma:\frac{1}{2}$ in terms of the other two.

\end{document}